\documentclass{article}

\usepackage{epsfig} 
\usepackage{amsmath} 
\usepackage{amssymb}  
\usepackage{amsthm}

\newtheorem{proposition}{Proposition}
\newtheorem{definition}{Definition}

\newtheorem{theorem}{Theorem}
\newtheorem{lemma}{Lemma}

\newtheorem{remark}{Remark}

\newcommand{\R}{{\mathbb{R}}}

\newcommand{\diag}{{\mathrm{diag}\,}}

\title{\LARGE \bf
Stability of continuous-time distributed consensus algorithms
}

\author{\Large Luc Moreau%
\thanks{The author was a Postdoctoral Fellow of the Fund for
        Scientific Research - Flanders (Belgium),
        associated to the SYsTeMS Research Group,
        Electrical energy, Systems and Automation Department,
        Ghent University, 9052 Ghent, Belgium.
        The author is currently employed at Sidmar, Ghent, Belgium.
        This paper presents research results of the Belgian Programme on
        Interuniversity Attraction Poles, initiated by the Belgian Federal
        Science Policy Office.  The scientific responsibility rests with its
        author.}
\thanks{Research performed while visiting the
        Department of Mechanical Engineering,
        Technical University of Eindhoven, 5600 MB Eindhoven, The Netherlands.}%
\\[2ex]
Ghent University\\ Belgium
}

\begin{document}

\maketitle

\begin{abstract}

We study the stability properties of
linear time-varying systems in continuous time
whose system matrix is Metzler with zero row sums. This class of systems
arises naturally in the context of distributed decision problems,
coordination and rendezvous tasks and synchronization problems.
The equilibrium set contains all states with identical state components.
We present sufficient conditions guaranteeing uniform exponential stability
of this equilibrium set, implying that all state components
converge to a common value as time grows unbounded.
Furthermore it is shown that this convergence result is robust with respect
to an arbitrary delay, provided that the delay affects only the off-diagonal
terms in the differential equation.\\[1ex]

\end{abstract}


\section{INTRODUCTION}

We consider a coordination problem for a network of
systems interacting via local coupling.
Each individual system in the network is assumed to have
simple integrator dynamics. The network coupling is allowed
to be time-dependent and non-bidirectional. We show that,
under very mild conditions, the individual state variables
of the systems in the network converge to a common value as
time grows unbounded.

Each system in the network may be a kinematic model
of a fully actuated vehicle and the corresponding state
variable may represent its position in one-dimensional space.
The coordination problem then asks to steer all vehicles to a common
rendezvous position---which is not {\it{a priori}} given but
will depend on the vehicles' initial positions.
The strategy studied here is of a distributed nature. No leader
or external coordinator is assumed to be present. Instead the
desired formation is formed as a consequence of the local interaction
between the individual vehicles.
Other vehicle models may be conceived. Each system in the network
could, for example, model a particle moving at constant speed
in two-dimensional space, its direction of motion being represented
by the system's state variable as in~\cite{ViCzJaCoSc:95,JaLiMo:03}.
More general models are possible, where the individual state
variables are not necessarily scalar but
vector-valued, representing the configuration of a vehicle
in three-dimensional space. Such higher dimensional extensions,
however, will not be considered in the present paper.
We refer
to~\cite{LeFi:01,SmHaLe:01,JuKr:techrep02,JuKr:04,SePaLe:03,LiPaPo:03,LiMoAn:03,LiBrFr:ieeetac04,LiFrMa:04,TaJaPa:04,NaLe:04,Ol:ieeetac04}
for the study of more general swarming models
and collective motion patterns.

In a different context of synchronization~\cite{St:03} each system
in the network could represent an oscillator and the individual
system's state
variable could determine the oscillation phase. In this case the
consensus problem corresponds to a synchronization problem for the
oscillators. The present model may then be interpreted as a linearized
version
of the more complex synchronization models studied, for example,
in~\cite{St:00,JaMoBa:04,AeRo:04}.

Finally, the individual systems' state variables could represent
an abstract decision variable. In this case the network equations
represent a distributed consensus algorithm for the purpose of
distributed decision making as in~\cite{OlMu:acc03,OlMu:cdc03,OlMu:ieeetac04,XiBo:scl04,ReBe:acc04}.
Distributed agreement problems have
a long history; see, for example, the book~\cite[section 4.6]{Se:81}
and the paper~\cite{De:74}.


The distributed consensus algorithm that we consider here has received considerable attention in the
recent literature~\cite{JaLiMo:03,LiBrFr:ieeetac04,OlMu:acc03,OlMu:cdc03,OlMu:ieeetac04,XiBo:scl04,ReBe:acc04,GaPa:03,ChWaMu:mtns04}.
Compared with these references,
a main contribution of the present paper is that we
do not impose restrictive assumptions on the network coupling topology (other than a
very mild connectivity assumption). We
allow for general non-bidirectional and time-dependent communication patterns.
Unidirectional communication is important in practical applications
and can easily be incorporated, for example, via broadcasting.  Also,
sensed information flow which plays a central role in schooling and
flocking is typically not bidirectional.  In addition, we do not
exclude loops in the coupling topology.  This means that,
typically, we are considering leaderless coordination rather than a
leader-follower approach.  (Leaderless coordination is also
considered, for example, in~\cite{FaMu:02a,FaMu:02b,OlMu:ifac02,BaLe:02,MoBaLe:lhmnlc03}).
Finally, we allow for time-dependent communication patterns which are
important if we want to take into account link failure and link
creation, reconfigurable networks and nearest neighbor coupling.

The present work may be seen as a continuation of the author's previous
work~\cite{Mo:cdc03,Mo:ieeetac04}. We extend the discrete-time results
of those papers towards continuous-time distributed algorithms.
Furthermore, we show that the continuous-time results presented here
are robust with respect to an
arbitrary communication delay, provided that the delay affects only those
variables that are actually being communicated between distinct systems
in the network.
Even though we restrict attention to linear systems in the present paper,
it is instructive to observe that our approach is of an inherently nonlinear nature,
with a non-quadratic Lyapunov function playing a central role.

The paper is organized as follows. Section~\ref{s:model} introduces the model.
Section~\ref{s:survey} gives an overview of various Lyapunov functions reported
in the literature. Sections~\ref{s:main} and~\ref{s:delay} present the main stability results,
respectively for the case without delay and the case with delay.

\section{MODEL}
\label{s:model}

We start with introducing the abstract
model that forms the subject of our study.
We consider
linear time-varying systems in continuous
time
\begin{equation}\label{e:linear}
\dot{x}=A(t)x
\end{equation}
where the system matrix~$A(t)$ is assumed to satisfy: (a) its off-diagonal
elements are
positive or zero and
(b) its row sums are zero.
Matrices with positive or zero off-diagonal
elements arise naturally in the study of
positive systems and Markov chains. They are
sometimes referred to as~{\em{Metzler}} matrices~\cite{Lu:79}
and their spectral properties are well documented
in the literature. The extra assumption concerning
zero row sums guarantees that each {\it{consensus}} state with
a common value for its state components corresponds to
an equilibrium state. What makes the present
study interesting is that the system matrix is
allowed to depend explicitly on time.

Before we continue our study of~\eqref{e:linear}
let us first provide some
background and motivation for considering
this class of systems; see also~\cite{JaLiMo:03,LiBrFr:ieeetac04,OlMu:acc03,OlMu:cdc03,OlMu:ieeetac04,XiBo:scl04,ReBe:acc04,GaPa:03,ChWaMu:mtns04}
where this model or variants of it have been studied.
We consider a finite network of systems.
(Throughout the paper both the individual entities making up the
network as well as the overall network are referred
to as `systems'. No confusion should arise.)
The state of each individual system in
the network, say system~$k$, is characterized by a real number~$x_k$.
The differential equations governing the
evolution of these variables depend on the network
interconnection.
Let us first consider a simple network structure
with only one coupling.
A link from
system~$l$ to system~$k$ contributes to the
dynamics as follows
$$
\dot{x}_k=a_{kl}(x_l-x_k)
$$
which makes the real variable~$x_k$
evolve towards~$x_l$ with a rate of change proportional
to the difference~$x_l-x_k$. The coefficient~$a_{kl}$ serves as
a measure for the coupling strength.
In matrix notation we have
\begin{equation*}\label{e:onelink}
\frac{\mathrm{d}}{\mathrm{d}t}
\begin{pmatrix}
\vdots\\
x_k\\
\vdots\\
x_l\\
\vdots
\end{pmatrix}
=
\begin{pmatrix}
&\vdots&&\vdots&\\
\dots&-a_{kl}&\dots&a_{kl}&\dots\\
&\vdots&&\vdots&\\
\dots&0&\dots&0&\dots\\
&\vdots&&\vdots&\\
\end{pmatrix}
\cdot
\begin{pmatrix}
\vdots\\
x_k\\
\vdots\\
x_l\\
\vdots
\end{pmatrix}.
\end{equation*}

Of course, we do not want to restrict attention to
single link networks.
When considering multiple
simultaneous couplings
we assume that their
effect is additive.
According to this additivity assumption the dynamics
of the network is described by
$\dot{x}=Ax$
where the system matrix has positive or zero
off-diagonal elements and zero row sums.
We may briefly refer to such matrices as
{\em{Metzler matrices with zero row sums}}.

In general, we want to allow the system matrix
to depend explicitly on time.
The system matrix may vary discontinuously with time, for
example, in the case of instantaneous
link creation or link
failure.
The system matrix may also vary continuously with time
when some coupling strengths increase or decrease
gradually. This could, for example, be the case when
the coupling strength depends on the distance between
systems in physical space.
Both continuous and discontinuous changes are
expected when the network of systems correspond to physical
entities interacting via nearest neighbor coupling.

\subsection{Graph representation}
\label{s:def}

It is illustrative to provide a graph
representation of the network coupling.
For this purpose we associate to any $n\times n$~Metzler
matrix with zero row sums a digraph (directed graph) as follows.
We introduce~$n$ nodes labelled consecutively from~$1$ to~$n$.
Each node represents a system in the network.
To each strictly positive off-diagonal
element of the matrix, say, on the
$k$th row and the
$l$th column
we associate an arc from node~$l$ towards
node~$k$ indicating that the state of system~$l$ influences
the dynamics of system~$k$. The orientation of the arc
corresponds to the direction of the information flow, namely from
system~$l$ to system~$k$.%
\footnote{Notice that we do not assign weights to
the arcs. It is straightforward
to construct a weighted digraph by taking
into account the precise
values of the matrix elements but this
is not needed in the present context.}

For some purposes we may wish to incorporate a threshold
when constructing a digraph from a given
Metzler matrix.
This is particularly useful when the system matrix depends
explicitly on time. In this case it may happen that some
matrix elements converge asymptotically to zero
and eventually become so small
that they hardly contribute to
the dynamics. It may be better to
neglect such small values when representing the network
coupling. This may be achieved by introducing a threshold
value~$\delta>0$. Only those
off-diagonal elements whose value is strictly
larger than~$\delta$
are then giving rise to an arc.
As this construction
will be instrumental in the statement of our main results
we provide the following formal definition (where we allow
$\delta=0$ in order to include the case of no threshold).
\begin{definition}[$\delta$-digraph associated to a matrix]
Consider an $n\times n$~Metzler matrix~$A$ with zero row sums.
The $\delta$-digraph ($\delta\geq 0$)
associated to~$A$ is the digraph with
node set~$\{1,\dots,n\}$ and with an arc from~$l$ to~$k$
($k\not =l$) if and only if the element of~$A$ on the
$k$th row and the
$l$th column is strictly larger than~$\delta$.
\end{definition}

It will turn out that the stability properties of~\eqref{e:linear}
may be related to the connectivity properties of
associated digraphs. In order to prepare for a precise formulation
of our results we introduce the following concept related to
connectivity. For a given digraph, we say that node~$l$ can be
reached from node~$k\not =l$ if there is a path from~$k$ to~$l$ in
the digraph which respects the orientation of the arcs.

\subsection{Stability and convergence}

We are interested in the convergence of the individual
state components to a common value. We may
reformulate this convergence property in terms of
classical stability concepts.
%
%
%
The equilibrium set of~\eqref{e:linear}
contains all {\it{consensus}} states with a common value for
their components. Convergence to a consensus
state may then be captured by requiring that the equilibrium set
of consensus states is attractive. More precisely we will be
focussing upon uniform exponential stability of the
set of consensus states. Sufficient conditions for uniform
exponential stability will be formulated in the remainder of
the paper.

If the system matrix does not depend explicitly on
time we may, of course, perform an eigenvalue
analysis.
A minor subtlety is introduced by the requirement that
the system matrix has zero row sums.
Indeed, it is straightforward to see that a matrix with
zero row sums has a zero
eigenvalue corresponding to
the right eigenvector~$(1,\dots,1)^\prime$
(where the superscript {\it{prime}} denotes transposition).
This trivial eigenvector determines the
equilibrium set~$\{x\in\R^n:x_1=\dots=x_n\}$.
The convergence property that we are aiming at
corresponds to all other eigenvalues having
strictly negative real part.
It is a well-known result implicitly contained
in~\cite[Sections~7.2 and~7.3]{Lu:79} that
all but one of the eigenvalues of a
Metzler matrix with zero row sums have
strictly negative real part, the only exception being the trivial
eigenvalue at zero, if and only if the
associated digraph has (at least) one node from which all other
nodes may be reached.%
\footnote{Notice that our construction of a digraph associated
to a Metzler matrix is different from the construction used
in the context of Markov chains. The difference concerns
the direction of the arcs, which is in the present paper
exactly the opposite as in the more familiar construction
used in the context of Markov chain theory. Stated in terms
of the latter convention we have that all but one of the eigenvalues of a
Metzler matrix with zero row sums have
strictly negative real part, the only exception being the trivial
eigenvalue at zero, if and only if the
associated digraph has (at least) one node which may be
reached from all other nodes.}
This provides a necessary and
sufficient condition for the convergence of the
individual systems' state variables to a common value.
In addition, the spectrum of the system matrix
provides a quantitative measure for the
speed of convergence; cfr.~\cite{TaJaPa:04,OlMu:acc03,OlMu:cdc03,XiBo:scl04}.
A further advantage of spectral analysis is that,
at least in principle,
it extends to the case of delayed communication,
although the explicit determination and characterization
of the spectrum may become more involved. (A useful result
in this context is the Gersgorin theorem~\cite{HoJo:88}; see, e.g., \cite{OlMu:cdc03,OlMu:ieeetac04}.)
An obvious disadvantage of this approach is that it is
limited to time-independent system matrices and that the
results are only locally valid when applied to nonlinear
systems.

\section{LYAPUNOV FUNCTIONS}
\label{s:survey}

Various approaches are available to analyse
the convergence properties of~\eqref{e:linear}.
In particular,
several Lyapunov functions have been proposed
in the literature.
The existence of a Lyapunov function may yield
useful insight in the system's qualitative behavior.
In addition, a Lyapunov function
may point towards generalizations of the present theory to
more complex coordination models, possibly involving
nonlinear and higher order dynamics. In this sense,
the present survey, although restricted to the particular system~\eqref{e:linear},
may be used to catalogue many
of the recent results on coordination and synchronization.

\subsection{Gradient flow}

If the system matrix is time-independent and symmetric,
$A(t)=A=A^\prime$ for all~$t$,
it is convenient to rewrite~\eqref{e:linear} as
\begin{equation}\label{e:gradient}
\dot{x}=-\nabla V(x)
\end{equation}
where we have introduced the potential
function~$V(x)=-x^\prime Ax/2$ and where
$\nabla V(x)$ denotes the gradient of~$V$ at~$x$:
$\nabla V(x)=(\partial V/\partial x_1,\dots,\partial V/\partial x_n)^\prime$.
It is clear that the equilibrium points
of~\eqref{e:gradient} correspond to stationary points
of~$V$ and that away from these points the
potential is strictly decreasing with time:
$$
\dot{V}(x)=-\|\nabla V(x)\|^2.
$$
This reformulation of the dynamics may serve
as a starting point of a stability analysis
which relates the stability properties of~\eqref{e:gradient}
with the qualitative properties of the potential
function~$V$.

In principle a gradient description is not restricted
to the case of time-independent system matrices,
but its usefulness is
limited when
the system matrix depends explicitly on time.
Indeed, in this case the potential function~$V(t,x)=-x^\prime A(t)x/2$
depends explicitly on time and
this gives rise to an extra term when evaluating
the rate of change of~$V$ with time:
$$
\dot{V}(t,x)=-\|\nabla V(t,x)\|^2-x^\prime \dot{A}(t)x/2.
$$
(Here we have assumed differentiability of $A(t)$.)
As a consequence the potential function is not longer
guaranteed to decrease along non-equilibrium
solutions, unless additional conditions are imposed.%
\footnote{One may, for example, impose the condition that
$x(t)\dot{A}(t)x(t)=0$ which may be interpreted as
a state-dependent
constraint on the time-evolution of the system matrix.}
An important advantage of this approach is that it extends to
nonlinear interactions, possibly yielding (almost) global convergence
results. In addition, this approach
extends to networks of second order dynamic models, where the potential
function should be augmented with kinetic energy terms in
order to obtain a non-increasing energy function; see~\cite{LeFi:01,TaJaPa:04,OlMu:ifac02}.
Finally we mention that, in the case of symmetric coupling, a stability
analysis based on {\it contraction analysis} is proposed in~\cite{SlWaRi:mtns04}.

\subsection{Sum of squares}

One may try to capture the convergence properties
of~\eqref{e:linear} by considering the
sum of squares
$S(x)=x_1^2+\dots+x_n^2$ as a candidate-Lyapunov function.
At first sight this may seem a surprising choice since
this function is not positive definite with respect to
the equilibrium set of all consensus states
but its usefulness is illustrated next.
The time-derivative of~$S(x)$ along the trajectories
of~\eqref{e:linear} is given
by~\(
\dot{S}(t,x)=x^\prime(A(t)+A^\prime(t))x
\)
and
thus we conclude that $S$~is non-increasing
with time
if and only if the symmetric matrix~$A(t)+A^\prime(t)$
is negative
semi-definite.
In this respect, the following result is very useful.
It states that, one, negative semi-definiteness of~$A(t)+A^\prime(t)$
is equivalent to~$A(t)$ having zero column sums. It also shows that
the sum of squares function is only a very simple representative
of a much broader class of Lyapunov functions that are all
non-increasing along the solutions of~\eqref{e:linear}.
The proof of this puzzling result is very short and elegant,
and is essentially contained in a 1976 paper by J.C.~Willems~\cite{Wi:76}.
\begin{proposition}\label{p:soq}
Let~$A$ be Metzler with zero row sums.
The following statements are equivalent:
\begin{enumerate}
\item\label{1}
The sum of squares function~$S(x)$ is non-increasing along the trajectories of~$\dot{x}=Ax$;
\item\label{2}
$A+A^\prime$ is negative semi-definite;
\item\label{3}
$A$ has zero column sums;
\item\label{4}
Every convex function~$V(x)$ invariant under coordinate permutations
is non-increasing along the trajectories of~$\dot{x}=Ax$.
\end{enumerate}
\end{proposition}
\begin{proof}
The fact that \ref{3} implies~\ref{4} follows from~\cite[Theorem 1(iii)]{Wi:76}.
We repeat the proof of that paper for reasons of completeness.
We show that~$V(\exp(At)x)\leq V(x)$ for every~$t\geq 0$.
It is not difficult to see that~$\exp(At)$ is a doubly stochastic
matrix; that is, a non-negative matrix with all row sums and column sums
equal to~$1$. According to a famous result by Birkhoff~\cite{HoJo:88}
a doubly stochastic matrix may be decomposed into a convex combination
of finitely many permutation matrices $\exp(At)=\sum_{i}\lambda_i P_i$ where
the coefficients~$\lambda_i$ are non-negative and sum up to one and
where the matrices~$P_i$ have exactly one entry equal to $1$ on each
row and column and all other entries equal to zero.
By convexity of~$V(x)$ and invariance under coordinate permutations
we have
\[
V(\exp(At)x)=V(\sum_{i}\lambda_i P_ix)\leq\sum_{i}\lambda_i V(P_ix)=\sum_{i}\lambda_i V(x)=V(x).
\]

The implications ``\ref{4} implies \ref{1} implies \ref{2}'' are trivial.
It remains to show that \ref{2} implies \ref{3}---a result which is
presumably well-known.
Let us denote by~$e$ the
vector~$(1,\dots,1)^\prime$. Using this notation we may
write the assumption that $A$~has zero row sums in
a very concise form, $Ae=0$. Similarly, the
requirement that $A$ has zero column sums corresponds
to~$e^\prime A=0$. In order to prove this, we
calculate~$x^\prime(A+A^\prime)x=2x^\prime Ax$ for the particular
case of~$x=e+\varepsilon A^\prime e$ where~$\varepsilon$ is
a small real number.
Using the row assumption~$Ae=0$ we have the following expression for~$x^\prime Ax$
\begin{equation*}
\begin{split}
(e^\prime+\varepsilon e^\prime A)A(e+\varepsilon A^\prime e)
&=
e^\prime A e+
\varepsilon ( e^\prime A A e + e^\prime A A^\prime e )
+ {0}(\varepsilon^2)\\
&=
\varepsilon e^\prime A A^\prime e
+ {0}(\varepsilon^2)\\
&=
\varepsilon |A^\prime e|^2
+ {0}(\varepsilon^2).
\end{split}
\end{equation*}
By requiring that this expression is negative or
zero for all~$\varepsilon$ we conclude that, necessarily,
$e^\prime A=0$ and thus that
the matrix~$A$ has zero column sums.
\end{proof}

We thus see that the sum of squares function or, more generally,
every convex function invariant under coordinate permutations is
non-increasing along the trajectories of~\eqref{e:linear}
if and only if the Metzler matrix~$A(t)$ has zero column
sums in addition to zero row sums.
This observation may serve as a starting
point for a Lyapunov/LaSalle stability analysis.

The sum of squares Lyapunov function has an appealing
interpretation in terms of passivity theory.
In order to clarify this point, we view the linear system
$\dot{x}=A(t)x$ as a feedback interconnection of a
multi-input multi-output system consisting of $n$~individual
integrators~$\dot{x}_k=u_k$ with a static map~$u=A(t)x$.
The integrators are
passive input-output systems with respective
storage functions~$x_i^2/2$. The feedback interconnection is
stable if (minus) the static map is also passive; that is, if
$x^\prime A(t)x\leq 0$ for all~$x$ and~$t$, which is precisely
the requirement we recognize from the above
discussion. In this case the
sum of the individual storage functions~$x_1^2+\dots +x_n^2$
(modulo a factor~$1/2$)
serves as a non-increasing Lyapunov function.
The present situation is only a very simple manifestation
of the interplay between passivity ideas and coordination
and synchronization tasks. We refer to the
paper~\cite{StSe:cdc03}
for more on this intimate relationship.

Apparently by coincidence, the
condition that the Metzler matrix~$A(t)$ has zero column sums
in addition to zero row sums implies that the sum~$x_1+\dots+x_n$
is constant along the solutions of~\eqref{e:linear}.
Indeed, the sum~$x_1+\dots+x_n$ may be written as~$e^\prime x$
and its time-derivative $(\mathrm{d}/\mathrm{d}t)e^\prime x=e^\prime A(t)x$
vanishes for all~$x$ when~$e^\prime A(t)=0$.
We may thus add any (nonlinear) function of
$x_1+\dots+x_n$ to $S(x)$
without changing the time-derivative. In particular, we may
replace the sum of squares Lyapunov function by
\begin{align*}
\tilde{S}(x)
&=\left(x_1-\frac{x_1+\dots+x_n}{n}\right)^2+\dots+\left(x_n-\frac{x_1+\dots+x_n}{n}\right)^2\\
&=S(x)-2(x_1+\dots+x_n)\frac{x_1+\dots+x_n}{n}+n\left(\frac{x_1+\dots+x_n}{n}\right)^2
\end{align*}
which satisfies $\dot{S}(t,x)=\dot{\tilde{S}}(t,x)$.
Unlike the original sum of squares function~$S(x)$,
the modified function~$\tilde{S}(x)$ has the appealing feature that it
is in fact positive definite with respect to the equilibrium set of consensus
states~$x_1=\dots=x_n$. This Lyapunov function is
used in~\cite{OlMu:cdc03,OlMu:ieeetac04} when distributed
consensus algorithms are studied that leave the average value of the state components invariant.

Finally, as has been recognized a long time ago, it is worth noticing
that the requirement of zero
row sums may be relaxed. In order to show this, let us start
with rewriting the zero column sum requirement in the
concise form~$e^\prime A=0$, as we did before.
The relaxation to be discussed here replaces the condition
$e^\prime A=0$ by the more general requirement~$p^\prime A=0$ where
$p$ is any vector with non-negative components.
It has been shown~\cite{Wi:76} that, if the system matrix
is Metzler with zero row sums and satisfies~$p^\prime A=0$
then every function of the
form~$p_1f(x_1)+\dots+p_nf(x_n)$ with $f$ an arbitrary convex function of
one real variable, is non-increasing along~$\dot{x}=Ax$.

To summarize, the Lyapunov functions introduced in this section enable
a stability analysis of distributed consensus algorithms. Advantages are
that the coupling
matrix is allowed to be time-dependent and that the approach may be extended
to incorporate nonlinear interactions, as in~\cite{JaMoBa:04}.
Furthermore, the passivity interpretation points towards far-reaching
generalizations involving higher-order models; see~\cite{StSe:cdc03}.
An important limitation, however, is that
restrictions need to be imposed on the network coupling (balance conditions in terms
of constraints on the column sums
of the coupling matrix). A simple leader-follower architecture
with an alternating
leader-follower relation as in~$\dot{x}=A(t)x$ with $A(t)$ switching
between
\[
\begin{pmatrix} -1 & 1 \\ 0 & 0\end{pmatrix}\mbox{ and }
\begin{pmatrix}  0 & 0 \\ 1 &-1\end{pmatrix},
\]
for example,
falls outside the scope of this method.

\subsection{Contraction property}

A third and last approach that we mention is based on
the {\it{contraction}} property%
\footnote{This terminology has been borrowed
from the text book~\cite{Se:81} and should not be confused with the
notion of contraction analysis advocated, for example, in~\cite{SlWaRi:mtns04,LoSl:98}.}
of~\eqref{e:linear}.
Recall that the dynamics of one component, say~$x_k$, is given by
\[
\dot{x}_k=\sum_{l=1,\,l\not = k}^n a_{kl}(t)(x_l-x_k)
\]
with~$a_{kl}(t)\geq 0$ for $k\not = l$.
This means that~$\dot{x}_k$ is a non-negative
linear
combination of the differences~$x_l-x_k$.
In other words,
each state component moves in the direction of the
other state components. It is thus clear that
$\max\{x_1,\dots,x_n\}$ is a non-increasing function of time.
Likewise $\min\{x_1,\dots,x_n\}$ cannot decrease.
We may combine these two properties by
introducing the Lyapunov function~$V(x)=\max\{x_1,\dots,x_n\}-\min\{x_1,\dots,x_n\}$
which is positive
definite with respect to the desired equilibrium
set~$\{x\in\R^n:x_1=\dots=x_n\}$ and
non-increasing along
the solutions of~\eqref{e:linear}.

The contraction property mentioned here has been known
for a long time, in particular in the context of discrete-time systems.
It serves as a starting point for the ergodicity analysis of
products of stochastic matrices~\cite{Se:81}.
As a matter of fact, even the application to
distributed consensus algorithms is already mentioned in that
reference; see~\cite[section~4.6]{Se:81}.


We illustrate below that, with the aid of the Lyapunov function~$V(x)=\max\{x_1,\dots,x_n\}-\min\{x_1,\dots,x_n\}$,
we are able to establish convergence under very mild assumptions concerning the network coupling.
Unlike the previous approaches, we do not need to assume symmetric coupling nor do we
need to require detailed balance conditions such as zero column sums.
The only assumption we have to make is a very mild connectivity assumption
over time intervals.
The obtained results complement and extend the discrete-time results
reported by the author in~\cite{Mo:cdc03,Mo:ieeetac04}.

\section{CONVERGENCE RESULT}
\label{s:main}

\begin{theorem}\label{t:main}
Consider the linear system
\begin{equation}\label{e:linearsys}
\dot{x}=A(t)x
\end{equation}
Assume that the system matrix is a bounded
and piecewise continuous function of time.
Assume that, for every time~$t$, the system
matrix is Metzler with zero row sums.
If there is an index~$k\in\{1,\dots,n\}$,
a threshold value~$\delta>0$ and an interval
length~$T>0$ such that
for all~$t\in\R$ the $\delta$-digraph
associated to
\begin{equation}\label{e:matrixint}
\int_{t}^{t+T}A(s)\,\mathrm{d}s
\end{equation}
has the property that all nodes
may be reached from the node~$k$, then
the equilibrium set of {\it{consensus}}
states is uniformly exponentially stable.
In particular,
all
components of any solution~$\zeta(t)$
of~\eqref{e:linearsys} converge
to a common value as~$t\rightarrow\infty$.
\end{theorem}
\begin{remark}
See section~\ref{s:def} for the precise
definition of a $\delta$-digraph and
a node which may be reached from another node.
\end{remark}

Before we give a proof we first provide
a brief discussion of Theorem~\ref{t:main}.
The node~$k$
may be interpreted as a {\it{dominant}}
system in the network in the sense that its
initial state value
will influence the final consensus value that will
be reached. This, however, does not necessarily mean that we
are dealing with a leader-follower
approach since there may be
many other dominant systems in the network
which will all contribute to the final consensus
value.

The connectivity condition featuring in
the statement of the theorem is a very mild
condition. For example, it does not require
the existence of a system
communicating directly
to all other systems in the network at any given time instant.
Instead, it allows the communication from one system
to another to be indirect, involving intermediate systems.
Also, the required communication does not need
to occur instantaneously but may be spread over time.

Clearly the piecewise continuity assumption
is merely a technical condition which
guarantees a proper meaning of~\eqref{e:linearsys}--\eqref{e:matrixint}.
In contrast to this, the boundedness assumption
plays an instrumental role in the proof
of Theorem~\ref{t:main}. Without this boundedness
assumption, the conclusion of the theorem
is no longer valid.

\subsection{Proof of Theorem~\ref{t:main}}
\label{s:proof}

We consider~$V(x)=\max\{x_1,\dots,x_n\}-\min\{x_1,\dots,x_n\}$
as a candidate Lyapunov function, where
$n$~is the dimension of the overall state space, $x\in\R^n$.
This function is positive
definite with respect to the desired equilibrium
set~$\{x\in\R^n:x_1=\dots=x_n\}$ and
non-increasing along
the solutions of~\eqref{e:linearsys}.
In general, the function~$V$ is not necessarily decreasing
at every time instant. However, it is possible to show
that~$V$ decreases over time intervals
of sufficient length.

In order to show that~$V$ decreases over time intervals of
sufficient length we calculate conservative but explicit
bounds for
the trajectories of~\eqref{e:linearsys}.
The following lemma
plays a crucial role in obtaining these
bounds.
\begin{lemma}\label{l:bound}
Consider the linear system
\begin{equation}
\dot{x}=A(t)x.\label{e:linearsys2}
\end{equation}
Assume that the system matrix is
a piecewise continuous function of time and that, for every time~$t$, the system
matrix is Metzler with zero row sums.
Let~$\zeta$ be a solution of~\eqref{e:linearsys2}
defined on the time interval~$[t_0,\,t_1]$.
Consider two nonempty disjoint sets~$G$ and~$H$
satisfying $G\cup H=\{1,\dots,n\}$ and
introduce the following numbers characterizing
the coupling relative with respect to~$G$ and~$H$
on the time interval~$[t_0,\,t_1]$
\begin{itemize}
\item
$a_{GH}=\sum_{k\in G,\,l\in H}\int_{t_0}^{t_1}a_{kl}(t)\,\mathrm{d}t$,
\item
$a_{HG}=\sum_{k\in H,\,l\in G}\int_{t_0}^{t_1}a_{kl}(t)\,\mathrm{d}t$,
\item
$a_{HH}=\sum_{k,l\in H,\,k\not =l}\int_{t_0}^{t_1}a_{kl}(t)\,\mathrm{d}t$,
\end{itemize}
where~$a_{kl}(t)$ denotes the element on the~$k$-th row and~$l$-th column
of the matrix~$A(t)$.
Introduce upper and lower
bounds $g^{\min}$, $g^{\max}$, $\mu^{\min}$ and $\mu^{\max}$ associated to the solution~$\zeta$ at time~$t_0$
\begin{itemize}
\item
$g^{\min}=\min_{k\in G}\{\zeta_k(t_0)\}$,
\item
$g^{\max}=\max_{k\in G}\{\zeta_k(t_0)\}$,
\item
$\mu^{\min}\leq\min_{k\in\{1,\dots,n\}}\{\zeta_k(t_0)\}$,
\item
$\mu^{\max}\geq\max_{k\in\{1,\dots,n\}}\{\zeta_k(t_0)\}$.
\end{itemize}
The following estimates hold for the solution~$\zeta$ at time~$t_1$:
\begin{itemize}
\item
every component~$\zeta_k(t_1)$, $k\in G$, is contained in the interval
\begin{equation}\label{e:intervalG}
[
\mu^{\min}+(g^{\min}-\mu^{\min}){\mathrm{e}}^{-a_{GH}},\,
\mu^{\max}-(\mu^{\max}-g^{\max}){\mathrm{e}}^{-a_{GH}}
],
\end{equation}
\item
at least one
component~$\zeta_k(t_1)$, $k\in H$, is contained in the interval
\begin{equation}\label{e:intervalH}
[
\mu^{\min}+(g^{\min}-\mu^{\min})\beta,\,
\mu^{\max}-(\mu^{\max}-g^{\max})\beta
],
\end{equation}
with
\begin{equation}\label{e:beta}
\beta=
\frac{{\mathrm{e}}^{-a_{GH}}a_{HG}/|H|}{1+{\mathrm{e}}^{a_{HH}}a_{HG}/|H|+{\mathrm{e}}^{a_{HH}}a_{HH}}.
\end{equation}
where~$|H|$ denotes the number of elements in~$H$.
\end{itemize}
\end{lemma}
The first estimate may be interpreted as a growth restriction
on the components~$\zeta_k(t_1)$, $k\in G$, and the second estimate
may be interpreted as a trapping region for at least one
of the components~$\zeta_k(t_1)$, $k\in H$. It is important to notice
that the estimates provided in~\eqref{e:intervalG}--\eqref{e:beta}
depend explicitly on the system matrix~$A(t)$ via the
numbers~$a_{GH}$, $a_{HG}$ and $a_{HH}$.
For each of these three numbers an upper bound
is available because the system matrix is assumed to be
a bounded function of time. In addition, for the second of these numbers,
a lower bound is available if appropriate connectivity conditions are
imposed, thus effectively yielding a lower bound for the number~$\beta$
featuring in~\eqref{e:intervalH}--\eqref{e:beta}.%
\footnote{It is left as an exercise to the reader to show
that upper bounds
on~$a_{GH}$ and~$a_{HH}$ and a lower bound on~$a_{HG}$ yield a lower
bound on~$\beta$.}

\begin{proof}[Proof of Lemma~\ref{l:bound}]
We start with proving the first of the two estimates; cfr.~\eqref{e:intervalG}.
Notice that, initially at time~$t_0$, every component~$\zeta_k$, $k\in G$,
is contained in the interval~$[
g^{\min},
g^{\max}
]$.
A component~$\zeta_k$, $k\in G$, can only leave this interval
when it is attracted by some of the other components~$\zeta_k$, $k\in H$.
These other components~$\zeta_k$, $k\in H$, (as a matter of fact, {\em{all}} components)
are themselves for all times~$t\geq t_0$ confined to the interval~$[
\mu^{\min},\,
\mu^{\max}
]$.
In view of this, it is not difficult to see that
the accumulative communication (attraction)~$a_{GH}$ from~$H$ to~$G$ yields a growth restriction
on the components~$\zeta_k$, $k\in G$.
Intuitively, it is clear that a maximal growth is realized if all communication
from~$H$ to~$G$ is actually coming from one or more components~$\zeta_k$, $k\in H$,
situated at one of the extremal positions, say, $\mu^{\max}$ and directed to
one component~$\zeta_k$, $k\in G$, originally situated at~$g^{\max}$.
In this situation, the component~$\zeta_k$, $k\in G$, receiving all communication
would end up at position~$\mu^{\max}-(\mu^{\max}-g^{\max}){\mathrm{e}}^{-a_{GH}}$.
This reasoning provides some intuition for the upper bound
in expression~\eqref{e:intervalG}.

A more formal argument goes as follows.
Let us denote~$\xi(t)=\max_{k\in G}\{\zeta_k(t)\}$ and let us
consider the time-evolution
of this function.
For almost every time instant~$t$ we have that~$\xi(t)=\zeta_\kappa(t)$
and~$\dot\xi(t)=\dot\zeta_\kappa(t)$ for some
(so-called extremal) index~$\kappa\in G$, possibly depending on~$t$.
The following holds for this extremal index
\begin{align*}
\dot\zeta_\kappa
&=\sum_{l\in G\setminus\{\kappa\}}a_{\kappa l}(t)(\zeta_l-\zeta_\kappa)+\sum_{l\in H}a_{\kappa l}(t)(\zeta_l-\zeta_\kappa)\\
&\leq\sum_{l\in H}a_{\kappa l}(t)(\mu^{\max}-\zeta_\kappa),
\end{align*}
where we used the fact that~$\zeta_l\leq\zeta_\kappa$ for all~$l\in G\setminus\{\kappa\}$
and~$\zeta_l\leq\mu^{\max}$ for all~$l\in H$.
We don't need to know the precise index~$\kappa$ since it always holds that
$\sum_{l\in H}a_{\kappa l}(t)\leq \sum_{k\in G,\,l\in H}a_{kl}(t)$.
We conclude that, at almost every time~$t$,
\[
\dot\xi\leq\sum_{k\in G,\,l\in H}a_{kl}(t)(\mu^{\max}-\xi).
\]
Hence, by integration, we have at time~$t_1$
\begin{align*}
\xi(t_1)&
\leq\mu^{\max}-(\mu^{\max}-\xi(t_0)){\mathrm{e}}^{-\sum_{k\in G,\,l\in H}\int_{t_0}^{t_1}a_{kl}(t)\,\mathrm{d}t}\\
&\leq\mu^{\max}-(\mu^{\max}-g^{\max}){\mathrm{e}}^{-a_{GH}}.
\end{align*}
This explains the upper bound in expression~\eqref{e:intervalG}.
The lower bound may, of course, be demonstrated similarly.

We now proceed with proving the second estimate of Lemma~\ref{l:bound}.
Suppose that at the final time~$t_1$ every~$\zeta_k$, $k\in H$, is situated
either in the interval~$[\mu^{\min},\,h^{\min})$ or in
the interval~$(h^{\max},\,\mu^{\max}]$, where~$h^{\min}$ and~$h^{\max}$
are two real numbers sufficiently close to~$\mu^{\min}$ and~$\mu^{\max}$, respectively.
By a similar reasoning as above (growth restriction) we conclude that,
initially at time~$t_0$
every~$\zeta_k$, $k\in H$, has to be situated
either in the interval
\begin{equation}\label{e:intervalHl}
[\mu^{\min},\,\mu^{\min}+(h^{\min}-\mu^{\min})\mathrm{e}^{a_{HH}})
\end{equation}
or in the interval
\begin{equation}\label{e:intervalHr}
(\mu^{\max}-(\mu^{\max}-h^{\max})\mathrm{e}^{a_{HH}},\,\mu^{\max}].
\end{equation}
Here, we assume implicitly that the growth can only be caused
by attraction from other components in $H$, but not from components in $G$.
This implicit assumption is justified if the two intervals~\eqref{e:intervalHl} and~\eqref{e:intervalHr} do not overlap with
the trapping region~\eqref{e:intervalG} for components in $G$.
This is the case if we restrict attention to values for~$h^{\min}$ and~$h^{\max}$
sufficiently close to~$\mu^{\min}$ and~$\mu^{\max}$, respectively.
The above intervals~\eqref{e:intervalHl} and~\eqref{e:intervalHr}
hold not only at the initial time~$t_0$, but also
{\it{a fortiori}} during the complete interval~$[t_0,\,t_1]$.

Now we turn attention to the communication from~$G$ to~$H$. Since there are $|H|$~components in~$H$, at least one of these
components~$\zeta_k$, $k\in H$,
receives a total amount of~$a_{HG}/|H|$ communication during the interval~$[t_0,\,t_1]$.
Assuming that this particular component, whose index we denote by~$\kappa$, is
in the right interval~$(h^{\max},\,\mu^{\max}]$ at time~$t_1$ and thus in the right
interval~\eqref{e:intervalHr} during~$[t_0,\,t_1]$
(left would be similar) we have
\begin{align*}
\dot\zeta_\kappa
&=\sum_{l\in G}a_{\kappa l}(t)(\zeta_l-\zeta_\kappa)+\sum_{l\in H\setminus\{\kappa\}}a_{\kappa l}(t)(\zeta_l-\zeta_\kappa)\\
&\leq
\sum_{l\in G}a_{\kappa l}(t)
(-(\mu^{\max}-g^{\max}){\mathrm{e}}^{-a_{GH}}+(\mu^{\max}-h^{\max})\mathrm{e}^{a_{HH}})\\
&\qquad
+\sum_{l\in H\setminus\{\kappa\}}a_{\kappa l}(t)(\mu^{\max}-h^{\max})\mathrm{e}^{a_{HH}}
\end{align*}
where we used the bound~\eqref{e:intervalG} for components in~$G$ and the bound~\eqref{e:intervalHr}
for component~$\zeta_\kappa$ as well as the upper bound~$\mu^{\max}$ for components in~$H$.
After integration we obtain
\begin{align}
\zeta_\kappa(t_1)
&\leq
\mu^{\max}
+(a_{HG}/|H|)
(-(\mu^{\max}-g^{\max}){\mathrm{e}}^{-a_{GH}}+(\mu^{\max}-h^{\max})\mathrm{e}^{a_{HH}})\notag\\
&\qquad
+a_{HH}(\mu^{\max}-h^{\max})\mathrm{e}^{a_{HH}}.\notag 
\end{align}
It is clear that the upper bound we just obtained for~$\zeta_\kappa(t_1)$
is getting tighter (monotonically decreasing) as~$h^{\max}$ increases towards~$\mu^{\max}$. At some
value for~$h^{\max}$ this upper bound actually yields
a contradiction with our starting assumption, namely that~$\zeta_\kappa(t_1)$ is
contained in the interval~$(h^{\max},\,\mu^{\max}]$. The critical value for~$h^{\max}$
at which this contradiction first occurs is easily obtained by solving an algebraic
equation and this yields the upper bound in expression~\eqref{e:intervalH}. The lower
bound may, of course, be obtained similarly. This concludes the proof of the technical lemma.
\end{proof}

Let us now show how Lemma~\ref{l:bound} may
be used in order to conclude that the Lyapunov
function~$V$ decreases over time intervals of
length~$(n-1)T$. Consider an arbitrary solution~$\zeta$
of~\eqref{e:linearsys} and an arbitrary initial
time~$t_0$.
Introduce the singleton~$G\in\{1,\dots,n\}$
corresponding to (one of) the node(s) from which all
other nodes may be reached for the
$\delta$-digraphs associated to~\eqref{e:matrixint}.
Let~$H$ be~$\{1,\dots,n\}\setminus G$ and
let~$\mu^{\min}=\min\{\zeta_1(t_0),\dots,\zeta_n(t_0)\}$
and~$\mu^{\max}=\max\{\zeta_1(t_0),\dots,\zeta_n(t_0)\}$---a
choice which is compatible with the data introduced in
Lemma~\ref{l:bound}. Notice that~$V(\zeta(t_0))=\mu^{\max}-\mu^{\min}$.
A first application of Lemma~\ref{l:bound} enables
us to conclude that, at time~$t_0+T$,
the component of~$\zeta(t_0+T)$ corresponding to the
singleton~$G$ is contained in a compact interval of
the form~\eqref{e:intervalG}, which is strictly contained
in~$[\mu^{\min},\,\mu^{\max}]$. At the same time,
it enables us to conclude that at least one
of the other state components~$\zeta_k(t_0+T)$, $k\in H$,
belongs to an interval of the form~\eqref{e:intervalH},
which is also strictly contained in~$[\mu^{\min},\,\mu^{\max}]$.%
\footnote{%
Regarding this last statement, it is important to emphasize
that the connectivity condition that we are assuming
to hold provides a lower bound $a_{HG}>\delta$ which
is instrumental for our conclusion that the interval
of the form~\eqref{e:intervalH}
is strictly contained in~$[\mu^{\min},\,\mu^{\max}]$.}
We may exclude (the index of) this particular component
from the set~$H$ and include it in the set~$G$ and invoke
Lemma~\ref{l:bound} a second time---with the same, original values
for~$\mu^{\min}$ and~$\mu^{\max}$ but with the
modified sets~$G$ and~$H$. Continuing along these
lines, we conclude after a repetitive application of
Lemma~\ref{l:bound} that, at time~$t_0+(n-1)T$, all
components of~$\zeta(t_0+(n-1)T)$ are contained in
a compact interval which is strictly contained in the
original interval~$[\mu^{\min},\,\mu^{\max}]$.
This shows that the Lyapunov function~$V$ has decreased
over the time interval~$[t_0,\,t_0+(n-1)T]$.
The estimates that are obtained above
by a repetitive application
of Lemma~\ref{l:bound} are conservative but explicit.
They depend on the
numbers~$a_{GH}$, $a_{HG}$ and $a_{HH}$
as well as on the initial differences~$g^{\min}-\mu^{\min}$
and~$\mu^{\max}-g^{\max}$.

The discussion of the previous paragraph
holds for every arbitrary solution~$\zeta$.
It establishes a decrease condition
$$
V(\zeta(t_0+(n-1)T)-V(\zeta(t_0))<0.
$$
In addition, the above discussion provides conservative
but explicit estimates for this decrease.
Based on these explicit estimates we may conclude the
existence of a continuous function~$\gamma(x)$ which
is positive definite with respect to the desired
equilibrium set~$\{x\in\R^n:x_1=\dots=x_n\}$
and which satisfies
$$
V(\zeta(t_0+(n-1)T)-V(\zeta(t_0))<-\gamma(\zeta(t_0))
$$
for any solution~$\zeta$ and any time~$t_0$.
Fairly standard Lyapunov arguments may now be invoked
in order to conclude that the desired equilibrium set
is uniformly exponentially stable%
\footnote{If the differential equations would be nonlinear
then extra conditions would have to be imposed on
the function~$\gamma$ in order to conclude
uniform exponential stability, but in the present
context this is not needed, since uniform asymptotic
stability and uniform exponential stability are equivalent
for linear systems.},
thus concluding the
proof of Theorem~\ref{t:main}.

\section{COMMUNICATION DELAY}
\label{s:delay}

It is well-known that, in general,
unmodelled delay effects in a feedback mechanism
may destabilize an otherwise stable system.
This destabilizing effect of delay has been
well documented in the literature.
In the present context delay effects may arise naturally,
for example, because of the finite transmission speed
due to the physical characteristics
of the medium transmitting the information
(e.g.\ acoustic wave communication between underwater vehicles).

As before, let us first consider a simple network structure
with one single coupling.
If we assume that the delay affects only the
variable that is actually being transmitted from one system
to another then
it makes sense to assume
that
a link from
system~$l$ to system~$k$ contributes to the
dynamics as follows
$$
\dot{x}_k(t)=a_{kl}(x_l(t-\tau)-x_k(t))
$$
which makes the real variable~$x_k$
evolve towards the delayed variable~$x_l$ with a rate of change proportional
to the difference~$x_l(t-\tau)-x_k(t)$.
More generally, we consider the
delay differential equation
$$
\dot{x}(t)=\diag(A(t))x(t)+(A(t)-\diag(A(t)))x(t-\tau)
$$
instead of~\eqref{e:linearsys}. Here~$\diag(A)$ is the obvious
notation for the diagonal matrix obtained from~$A$ by setting
all off-diagonal entries equal to zero.
Observe that the delay~$\tau$ is only featuring in those terms
that correspond to the off-diagonal elements of~$A(t)$.
Here and in the
remainder of the paper $\tau$~is a fixed positive real number.
The stability result
of the previous section is robust with respect to an arbitrary delay~$\tau$.
\begin{theorem}\label{t:delay}
Consider the linear system
\begin{equation}\label{e:linearsysdelay}
\dot{x}(t)=\diag(A(t))x(t)+(A(t)-\diag(A(t)))x(t-\tau)
\end{equation}
with~$\tau>0$.
Assume that the system matrix~$A(t)$ is a bounded
and piecewise continuous function of time.
Assume that, for every time~$t$, the system
matrix is Metzler with zero row sums.
If there is~$k\in\{1,\dots,n\}$, $\delta>0$ and~$T>0$ such that
for all~$t\in\R$ the $\delta$-digraph
associated to
\begin{equation}\label{e:matrixintdelay}
\int_{t}^{t+T}A(s)\,\mathrm{d}s
\end{equation}
has the property that all nodes
may be reached from the node~$k$, then
the equilibrium set of {\it{consensus}}
states is uniformly exponentially stable.
In particular,
all
components of any solution~$\zeta(t)$
of~\eqref{e:linearsysdelay} converge
to a common value as~$t\rightarrow\infty$.
\end{theorem}

This result can actually be proven following the lines
of the proof of Theorem~\ref{t:main}.
The main ideas used in that
proof essentially remain valid even if a delay is featuring
in the equations. The most important modification is that now,
instead of the Lyapunov function~$V(x(t))=\max\{x_1(t),\dots,x_n(t)\}-\min\{x_1(t),\dots,x_n(t)\}$,
we consider the Lyapunov {\em{functional}}
\[
{\mathcal{V}}(x_t)=\max_{\sigma\in[t-\tau,\,t]}\{x_1(\sigma),\dots,x_n(\sigma)\}-\min_{\sigma\in[t-\tau,\,t]}\{x_1(\sigma),\dots,x_n(\sigma)\},
\]
where~$x_t$ is a standard notation from the theory of delay differential equations and denotes
the time function~$x(\sigma)$ on the interval~$\sigma\in[t-\tau,\,t]$.
Once this modification is in place, all the other modifications that have to
be made are straightforward and left to the reader.

\begin{remark}
Distributed consensus algorithms with delay have also been
considered in~\cite{OlMu:cdc03,OlMu:ieeetac04}. In those references
the delay is assumed to affect the off-diagonal terms as
well as the diagonal terms in the differential equation.
In that case a large delay value may, in fact, destabilize the algorithm.
\end{remark}

\section{Conclusion}

An important trend in systems theory is to view a complicated system as an
interconnection of simpler subsystems. Doing this, one hopes to be able
to relate the behavior of the overall system to the dynamics of the
simpler subsystems as well as the interconnection topology.
In the present study, we encountered a very special instance of this paradigm,
were each individual system in the network has trivial integrator dynamics.
We have shown that in this case, only very mild conditions need to be imposed
on the interconnection topology in order to draw relevant conclusions
for the overall system in terms of consensus reaching. These very mild conditions
allow the communication to be time-dependent and non-bidirectional. In addition,
arbitrary time-delays are also allowed. It would be of interest to identify
other classes of systems beyond the simple integrators considered here, where similar
mild conditions on the coupling topology suffice to draw conclusions about
the overall system behavior.

A possible contribution in this direction could come from the theory of
monotone systems~\cite{Sm:95,AnSo:03}. The model studied in the present paper is,
in fact, a monotone, cooperative system---an
observation which we have not explicitly exploited in the present paper.
It may be interesting to investigate how the present results relate to the
recent work on cooperative systems, monotone systems and small gain
theorems~\cite{AnSo:03,So:02}.
This may possibly constitute a fruitful line of research.



\addtolength{\textheight}{-3cm}   


\begin{thebibliography}{10}
\providecommand{\url}[1]{#1}
\csname url@rmstyle\endcsname
\providecommand{\newblock}{\relax}
\providecommand{\bibinfo}[2]{#2}
\providecommand\BIBentrySTDinterwordspacing{\spaceskip=0pt\relax}
\providecommand\BIBentryALTinterwordstretchfactor{4}
\providecommand\BIBentryALTinterwordspacing{\spaceskip=\fontdimen2\font plus
\BIBentryALTinterwordstretchfactor\fontdimen3\font minus
  \fontdimen4\font\relax}
\providecommand\BIBforeignlanguage[2]{{%
\expandafter\ifx\csname l@#1\endcsname\relax
\typeout{** WARNING: IEEEtran.bst: No hyphenation pattern has been}%
\typeout{** loaded for the language `#1'. Using the pattern for}%
\typeout{** the default language instead.}%
\else
\language=\csname l@#1\endcsname
\fi
#2}}

\bibitem{ViCzJaCoSc:95}
T.~Vicsek, A.~Czir\'ok, E.~Ben-Jacob, I.~Cohen, and O.~Schochet, ``Novel type
  of phase transitions in a system of self-driven particles,'' \emph{Phys.\
  Rev.\ Lett.}, vol.~75, no.~6, pp. 1226--1229, Aug. 1995.

\bibitem{JaLiMo:03}
A.~Jadbabaie, J.~Lin, and A.~S. Morse, ``Coordination of groups of mobile
  autonomous agents using nearest neighbor rules,'' \emph{{IEEE} Trans.\
  Automat.\ Control}, vol.~48, no.~6, pp. 988--1001, June 2003.

\bibitem{LeFi:01}
N.~E. Leonard and E.~Fiorelli, ``Virtual leaders, artificial potentials and
  coordinated control of groups,'' in \emph{Proceedings of the 40th IEEE
  Conference on Decision and Control}, 2001, pp. 2968--2973, {O}rlando,
  Florida, USA.

\bibitem{SmHaLe:01}
T.~R. Smith, H.~Hanssmann, and N.~E. Leonard, ``Orientation control of multiple
  underwater vehicles with symmetry-breaking potentials,'' in \emph{Proceedings
  of the 40th {IEEE} Conference on Decision and Control}, 2001, pp. 4598--4603,
  {O}rlando, {F}lorida, {USA}, {D}ecember~4--7.

\bibitem{JuKr:techrep02}
E.~Justh and P.~Krishnaprasad, ``A simple control law for {UAV} formation
  flying,'' 2002, {I}nstitute for Systems Research Technical Report TR 2002-38.

\bibitem{JuKr:04}
------, ``Equilibria and steering laws for planar formations,'' \emph{Systems
  Control Lett.}, vol.~52, pp. 25--38, 2004.

\bibitem{SePaLe:03}
R.~Sepulchre, D.~Paley, and N.~Leonard, ``Collective motion and oscillator
  synchronization,'' \emph{Proceedings of the Block Island Workshop on
  Cooperative Control}, June 2003.

\bibitem{LiPaPo:03}
Y.~Liu, K.~M. Passino, and M.~Polycarpou, ``Stability analysis of
  one-dimensional asynchronous swarms,'' \emph{{IEEE} Trans.\ Automat.\
  Control}, vol.~48, no.~10, pp. 1848--1854, Oct. 2003.

\bibitem{LiMoAn:03}
J.~Lin, A.~S. Morse, and B.~D.~O. Anderson, ``The multi-agent rendezvous
  problem,'' in \emph{Proceedings of the 42nd IEEE Conference on Decision and
  Control}, 2003, pp. 1508--1513, {M}aui, Hawaii, USA, December~9--12.

\bibitem{LiBrFr:ieeetac04}
Z.~Lin, M.~Broucke, and B.~Francis, ``Local control strategies for groups of
  mobile autonomous agents,'' \emph{{IEEE} Trans.\ Automat.\ Control}, vol.~49,
  no.~4, pp. 622--629, Apr. 2004.

\bibitem{LiFrMa:04}
Z.~Lin, B.~Francis, and M.~Maggiore, ``Necessary and sufficient graphical
  conditions for formation control of unicycles,'' \emph{{IEEE} Trans.\
  Automat.\ Control}, 2004, to appear.

\bibitem{TaJaPa:04}
H.~Tanner, A.~Jadbabaie, and G.~J. Pappas, ``Flocking in fixed and switching
  networks,'' \emph{Submitted for publication}, 2004, see also the Proceedings
  of the 42nd IEEE Conference on Decision and Control, Maui, Hawaii, USA,
  December~9--12, 2003.

\bibitem{NaLe:04}
S.~Nair and N.~E. Leonard, ``Stabilization of a coordinated network of rotating
  rigid bodies,'' in \emph{Proceedings of the 43rd {IEEE} Conference on
  Decision and Control}, 2004, submitted for publication.

\bibitem{Ol:ieeetac04}
R.~Olfati-Saber, ``Flocking for multi-agent dynamic systems: {A}lgorithms and
  theory,'' \emph{{IEEE} Trans.\ Automat.\ Control}, 2004, submitted for
  publication.

\bibitem{St:03}
S.~H. Strogatz, \emph{Sync: The Emerging Science of Spontaneous Order}.\hskip
  1em plus 0.5em minus 0.4em\relax Hyperion, 2003.

\bibitem{St:00}
------, ``From {K}uramoto to {C}rawford: exploring the onset of synchronization
  in populations of coupled oscillators,'' \emph{Physica {D}}, vol. 143, pp.
  1--20, 2000.

\bibitem{JaMoBa:04}
A.~Jadbabaie, N.~Motee, and M.~Barahona, ``On the stability of the {K}uramoto
  model of coupled nonlinear oscillators,'' in \emph{Proceedings of the
  American Control Conference}, 2004, submitted for publication.

\bibitem{AeRo:04}
D.~Aeyels and J.~Rogge, ``Stability of phase locking and existence of
  entrainment in networks of globally coupled oscillators,'' in
  \emph{Proceedings of the 6th IFAC symposium on Nonlinear Control Systems},
  Sept. 2004, to appear.

\bibitem{OlMu:acc03}
R.~Olfati-Saber and R.~M. Murray, ``Consensus protocols for networks of dynamic
  agents,'' in \emph{Proceedings of the American Control Conference}, June
  2003, pp. 951--956, {D}enver, Colorado, {USA}.

\bibitem{OlMu:cdc03}
------, ``Agreement in networks with directed graphs and switching topology,''
  in \emph{Proceedings of the {IEEE} Conference on Decision and Control}, Dec.
  2003, pp. 4126--4132, {M}aui, Hawaii, {USA}.

\bibitem{OlMu:ieeetac04}
------, ``Consensus problems in networks of agents with switching topology and
  time-delays,'' \emph{{IEEE} Trans.\ Automat.\ Control}, 2004, to appear.

\bibitem{XiBo:scl04}
L.~Xiao and S.~Boyd, ``Fast linear iterations for distributed averaging,''
  \emph{Systems Control Lett.}, vol.~53, pp. 65--78, 2004.

\bibitem{ReBe:acc04}
W.~Ren and R.~W. Beard, ``Consensus of information under dynamically changing
  interaction topologies,'' in \emph{Proceedings of the American Control
  Conference}, 2004, {B}oston, Massachusetts, USA.

\bibitem{Se:81}
E.~Seneta, \emph{Non-negative Matrices and {M}arkov Chains}, 2nd~ed., ser.
  Springer Series in Statistics.\hskip 1em plus 0.5em minus 0.4em\relax New
  York: Springer, 1981.

\bibitem{De:74}
M.~H. DeGroot, ``Reaching a consensus,'' \emph{Journal of the American
  Statistical Association}, vol.~69, no. 345, pp. 118--121, Mar. 1974.

\bibitem{GaPa:03}
V.~Gazi and K.~M. Passino, ``Stability analysis of swarms,'' \emph{{IEEE}
  Trans.\ Automat.\ Control}, vol.~48, no.~4, pp. 692--697, Apr. 2003.

\bibitem{ChWaMu:mtns04}
T.~Chu, L.~Wang, and S.~Mu, ``Collective behavior analysis of an anisotropic
  swarm model,'' in \emph{Proceedings of the 16th International Symposium on
  Mathematical Theory of Networks and Systems}, 2004, {L}euven, {B}elgium,
  {J}uly~5--9.

\bibitem{FaMu:02a}
J.~A. Fax and R.~M. Murray, ``Information flow and cooperative control of
  vehicle formations,'' in \emph{Proceedings of the {IFAC} World Congress},
  July 2002, {B}arcelona, Spain.

\bibitem{FaMu:02b}
------, ``Graph {L}aplacians and stabilization of vehicle formations,'' in
  \emph{Proceedings of the {IFAC} World Congress}, July 2002, {B}arcelona,
  Spain.

\bibitem{OlMu:ifac02}
R.~Olfati-Saber and R.~M. Murray, ``Distributed cooperative control of multiple
  vehicle formations using structural potential functions,'' in
  \emph{Proceedings of the {IFAC} World Congress}, July 2002, {B}arcelona,
  Spain.

\bibitem{BaLe:02}
R.~Bachmayer and N.~E. Leonard, ``Vehicle networks for gradient descent in a
  sampled environment,'' in \emph{Proceedings of the 41th Conference on
  Decision and Control}, 2002, pp. 112--117, {L}as Vegas, Nevada, USA,
  December~10--13.

\bibitem{MoBaLe:lhmnlc03}
L.~Moreau, R.~Bachmayer, and N.~E. Leonard, ``Coordinated gradient descent: A
  case study of {L}agrangian dynamics with projected gradient information,'' in
  \emph{Preprints of the 2nd {IFAC} workshop on {L}agrangian and {H}amiltonian
  methods for nonlinear control}, A.~Astolfi, A.~van~der Schaft, and
  F.~Gordillo, Eds., 2003, pp. 67--72, {S}eville, {S}pain, April~3--5.

\bibitem{Mo:cdc03}
L.~Moreau, ``Leaderless coordination via bidirectional and unidirectional
  time-dependent communication,'' in \emph{Proceedings of the 42nd {IEEE}
  Conference on Decision and Control}, 2003, pp. 3070--3075, {M}aui, Hawaii,
  {USA}, {D}ecember~9--12.

\bibitem{Mo:ieeetac04}
------, ``Stability of multi-agent systems with time-dependent communication
  links,'' \emph{{IEEE} Trans.\ Automat.\ Control}, 2004, accepted for
  publication.

\bibitem{Lu:79}
D.~G. Luenberger, \emph{Introduction to Dynamic Systems: Theory, Models and
  Applications}.\hskip 1em plus 0.5em minus 0.4em\relax John Wiley \& Sons,
  1979.

\bibitem{HoJo:88}
R.~A. Horn and C.~R. Johnson, \emph{Matrix Analysis}.\hskip 1em plus 0.5em
  minus 0.4em\relax Cambridge University Press, 1988, first published 1985.

\bibitem{SlWaRi:mtns04}
J.-J.~E. Slotine, W.~Wang, and K.~E. Rifai, ``Contraction analysis of
  synchronization in networks of nonlinearly coupled oscillators,'' in
  \emph{Proceedings of the 16th International Symposium on Mathematical Theory
  of Networks and Systems}, 2004, {L}euven, {B}elgium, {J}uly~5--9.

\bibitem{Wi:76}
J.~C. Willems, ``Lyapunov functions for diagonally dominant systems,''
  \emph{Automatica J. {IFAC}}, vol.~12, pp. 519--523, 1976.

\bibitem{StSe:cdc03}
G.-B. Stan and R.~Sepulchre, ``Dissipativity characterization of a class of
  oscillators and networks of oscillators,'' in \emph{Proceedings of the 42nd
  {IEEE} Conference on Decision and Control}, Dec. 2003, pp. 4169--4173,
  {M}aui, Hawaii, {USA}.

\bibitem{LoSl:98}
W.~Lohmiller and J.~Slotine, ``On contraction analysis for non-linear
  systems,'' \emph{Automatica J. {IFAC}}, vol.~34, no.~6, pp. 683--696, 1998.

\bibitem{Sm:95}
H.~Smith, \emph{Monotone dynamical systems: an introduction to the theory of
  competitive and cooperative systems}, ser. Mathematical Surveys and
  Monographs.\hskip 1em plus 0.5em minus 0.4em\relax {P}rovidence, {RI}:
  American Mathematical Society, 1995, vol.~41.

\bibitem{AnSo:03}
D.~Angeli and E.~D. Sontag, ``Monotone control systems,'' \emph{{IEEE} Trans.\
  Automat.\ Control}, vol.~48, no.~10, pp. 1684--1698, Oct. 2003.

\bibitem{So:02}
E.~D. Sontag, ``Asymptotic amplitudes and cauchy gains: a small-gain principle
  and an application to inhibitory biological feedback,'' \emph{Systems Control
  Lett.}, vol.~47, pp. 167--179, 2002.

\end{thebibliography}

\end{document}